\documentclass[ao]{iosart2x}

\usepackage{dcolumn}
\usepackage{graphicx}
\usepackage{amsfonts}
\usepackage{url}
\newcolumntype{d}[1]{D{.}{.}{#1}}

\firstpage{1}
\lastpage{19}
\volume{}
\pubyear{0}

\begin{document}
\begin{frontmatter} 

%
\title{The grounding for Continuum}

\runningtitle{Continuum}

\author[A]{\inits{S.}\fnms{Stanislaw} \snm{Ambroszkiewicz}\ead[label=e1]{sambrosz@gmail.com}%
\thanks{Corresponding author. \printead{e1}.}%
},
\runningauthor{S. Ambroszkiewicz}
\address[A]{Institute of Computer Science, \institution{Polish Academy of Sciences}, al. Jana Kazimierza 5, 01-248 Warsaw, \cny{Poland}\printead[presep={\\}]{e1}}

\begin{abstract}
It is a ubiquitous opinion among mathematicians that a real number is just a point in the line. If this rough definition is not enough, then a mathematician may provide a formal definition of the real numbers in the set theoretic and axiomatic fashion, i.e. via Cauchy sequences or Dedekind cuts, or as the collection of axioms  characterizing exactly (up to isomorphism) the set of real numbers as the complete and totally ordered Archimedean field.  
Actually, the above notions of the real numbers are abstract and do not have a constructive grounding. Definition of Cauchy sequences, and equivalence classes of these sequences explicitly use the actual infinity. The same is for Dedekind cuts, where the set of rational numbers is used as actual infinity. 
Although there is no direct constructive grounding for these  abstract notions, there are so called intuitions on which they are based. 
A rigorous approach to express these very intuition in a constructive way is proposed. It is based on the concept of the adjacency relation that seems to be a missing primitive concept in type theory. The approach corresponds to the intuitionistic view of Continuum proposed by Brouwer. The famous and controversial Brouwer Continuity Theorem is discussed on the basis of different principle than the Axiom of Continuity. 
\end{abstract}


\end{frontmatter}

\section{Introduction} 

In the XIX century and at the beginning of the XX century there was a common view that Continuum cannot be reduced to numbers, that is, Continuum cannot be identified with the set of the real numbers, and in general with a compact connected metric space. Real numbers are defined on the basis of rational numbers as equivalence classes of Cauchy sequences, or Dedekind cuts. 

The following citations support this view.
\begin{itemize}
\item David \cite{Hilbert1935}: {\em the geometric continuum is a concept in its own right and independent of number}. 
\item Emile Borel (in \cite{Troelstra}): {\em ...  had to accept the continuum as a primitive concept, not reducible to an arithmetical theory of the continuum [numbers as points, continuum as a set of points]}. 
\item Luitzen E. J. \cite{Brouwer1913} (and in \cite{Troelstra}): {\em The continuum as a whole was intuitively given to us; a construction of the continuum, an act which would create all its parts as individualized by the mathematical intuition is unthinkable and impossible. The mathematical intuition is not capable of creating other than countable quantities in an individualized way. [...] the natural numbers and the continuum as two aspects of a single intuition (the primeval intuition).}
\end{itemize}

In the mid-1950s  there were some attempts to comprehend the intuitive notion of Continuum by giving it strictly computational and constructive sense, i.e. by considering computable real numbers and computable functions on those numbers, see Andrzej  \cite{Grzegorczyk1955a} \cite{Grzegorczyk1955b} \cite{Grzegorczyk1957}, and Daniel  \cite{Lacombe1955}. These approaches were mainly logical and did not find a ubiquitous interest in Mathematics. 

If the concept of Continuum is different than the concept of number, then the problem of  reconstructing the computational grounding for the Continuum is important. 
For a comprehensive review and  discussion with a historical background, see \cite{feferman2009conceptions, bell2005divergent, Longo1999-LONTMC}. 

 Recently, see HoTT \cite{HoTT}, a type theory was introduced to homotopy theory in order to add computational and constructive aspects. However, it is based on Per Martin L\"{o}f's type theory that is a formal theory invented to provide intuitionist foundations for Mathematics.  The authors of HoTT  admit that there is still no computational grounding for HoTT. 

\cite{RobertHarper} :
{\em  ``...
And yet, for all of its promise, what HoTT currently lacks is a computational interpretation!  What, exactly, does it mean to compute with higher-dimensional objects? ...
 type theory is and always has been a theory of computation on which the entire edifice of mathematics ought to be built.  ... ``}
  
The Continuum  is defined  in HoTT as the real numbers via Cauchy sequences. 

Since the intuitive notion of Continuum is common for all humans (not only for mathematicians), the computational grounding of the Continuum (as a primitive type) must be simple and obvious. 

According to the American Heritage\textregistered Dictionary of the English Language: 
{\em ``Continuum is a continuous extent, succession, or whole, no part of which can be distinguished from neighboring parts except by arbitrary division''}.
 
Intuitively, a Continuum can be divided into finitely many parts. Each of the parts is the same Continuum modulo size. Two adjacent parts can be united, and the result is a Continuum. It is not surprising that this very common intuition of the Continuum goes back to Aristotle (384--322 BCE) and his conception of Continuum. Recently, the Aristotelian Continuum has got a remarkable research interest focused on axiomatic characterization of Continuum related to real line (the set of real numbers), and two dimensional Euclidean space  $\mathbb{R}^2$, see  \cite{roeper2006}, 
\cite{hellman2013classical}, 
\cite{linnebo2015aristotelian}, 
\cite{doi:10.1093/philmat/nkv021}. Roughly, the axiomatic theories define region-based (pointless) topologies. The proposed axiomatic systems are categorical. This means that each of them describes the unique model up to isomorphism.  
The theories make essential use of actual infinity, so that these approaches are not constructive. There is also an axiomatic approach that does not use actual infinity, see  \cite{linnebo2017actual}, however, it is limited. 

The grounding of Continuum presented in the next sections may be seen as constructive models of these formal theories. 
It is based on the intuitionistic  notion of Continuum proposed by Brouwer. 


The grounding is extended to the generic notion of Continuum and various interesting related topological spaces. The main idea is simple and is based on W-types augmented with simply patterns of adjacency relation. Finally, the famous and controversial Brouwer Continuity Theorem is proved on the basis of a different principle than the Axiom of Continuity  originally used to prove the Theorem by Brouwer.

%
%
\section{Examples } 
\label{examples}

Since the closed unit interval $[0;\ 1]$, a subset of the set of real numbers  $\mathbb{R}$, is a mathematical example of Continuum, let us follow this interpretation. Let the $n$-th dimensional Continuum be interpreted as the $n$-th dimensional unit cube $[0;\ 1]^n$, that is, the Cartesian product of $n$-copies of the unit interval. Note that 3-dimensional cube ($3$-cube, for short) is the only one regular polyhedron filling 3D-space. 

Let us fix the dimension $n$ and consider a unit cube. The cube may be divided into parts (smaller $n$-cubes) in many ways. However, the most uniform division is to divide it into $2^n$ the same parts, where $n$ is the dimension of the cube. 
Each of the parts may be divided again into $2^n$ the same sub-parts (even smaller $n$-cubes), and so on, finitely many times. There is a natural relation of adjacency between the elementary cubes. Two cubes (parts) are adjacent if their intersection (as sets) is a cube of dimension $n-1$. 

Our purpose is to liberate (step by step) our reasoning from the Euclidean interpretation of Continuum as $n$-cubes.
It seems that only the following properties are essential. 
\begin{itemize}
\item Any part is {\em the same} (modulo size) as the original unit cube. It can be divided arbitrary many times into smaller parts.  The same patter of partition, and the same adjacency pattern between the resulting parts is repeated in the consecutive divisions. 
\item Any part (except the border parts) has the same number of its adjacent parts and the same adjacency structure. 
\end{itemize}
Hence, only the following aspects are important:
\begin{itemize}
\item  the adjacency relation between parts resulting from consecutive divisions; 
\item  fractal-like structure of the unit cube, i.e. its  parts are similar (isomorphic) to the cube; 
\item the partition pattern is the same for all divisions;
\item the pattern of adjacency between sub-parts, after consecutive divisions, is one and the same.   
\end{itemize}
\begin{figure}[h]
	\centering
	\includegraphics[width=0.6\textwidth]{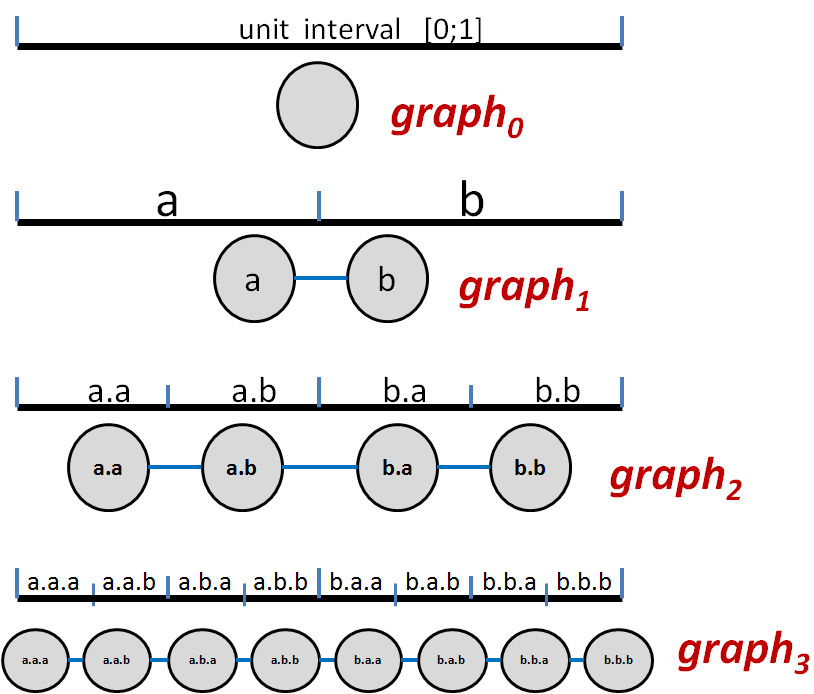}
	\caption{Consecutive partitions of Continuum of dimension $1$, i.e. the unit interval $[0; 1]$, and the corresponding graphs as adjacency relations between parts}
	\label{unit-interval}
\end{figure} 
\begin{figure}[h]
	\centering
	\includegraphics[width=0.6\textwidth]{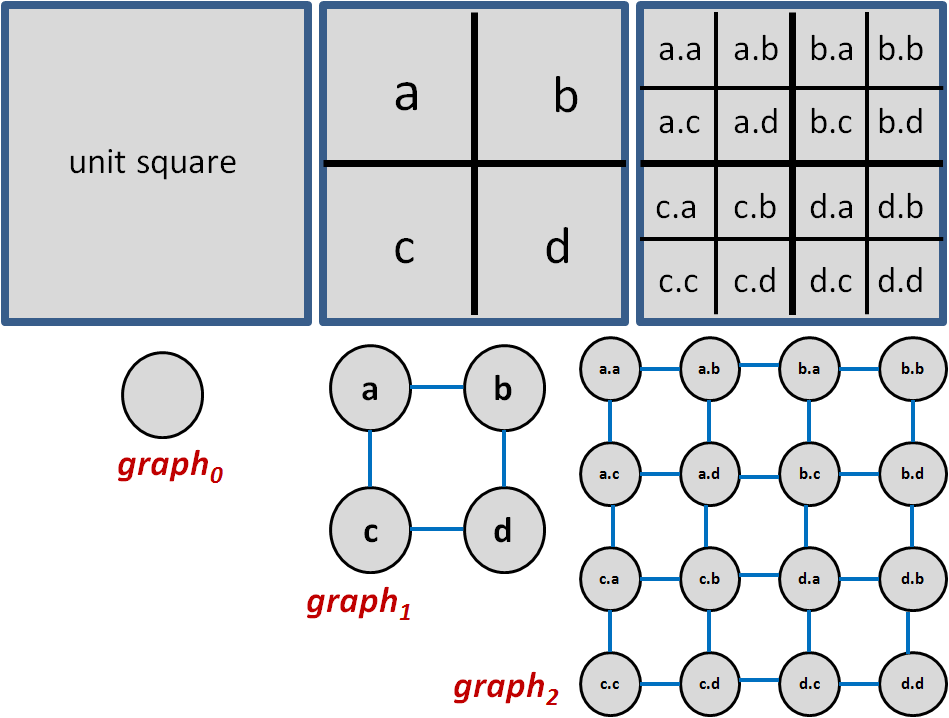}
	\caption{Consecutive partitions of Continuum of dimension $2$, i.e. the unit square $[0; 1]^2$, and the corresponding graphs as adjacency relations between parts}
	\label{unit-square}
\end{figure} 
\begin{figure}[h]
	\centering
	\includegraphics[width=0.6\textwidth]{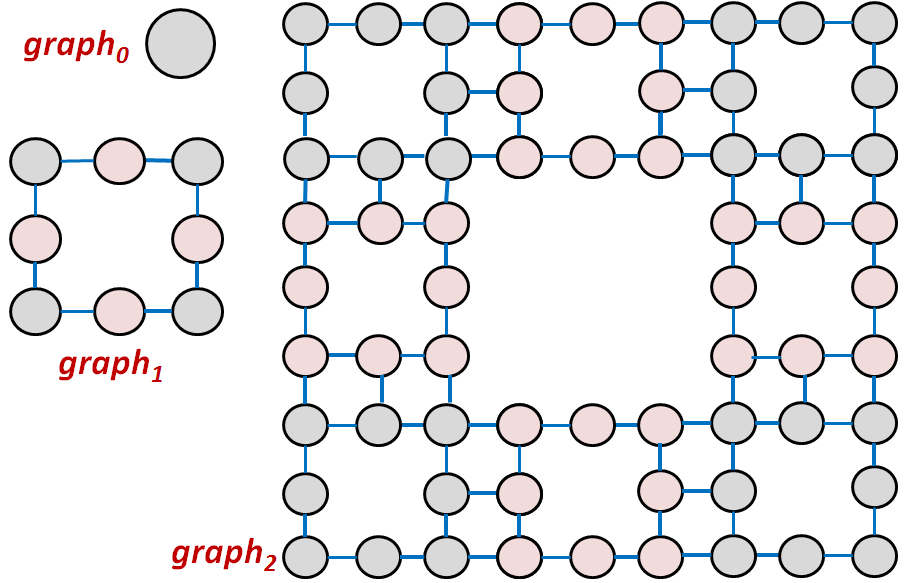}
	\caption{Consecutive partitions of an abstract Continuum corresponding to the Sierpinski carpet (a fractal), and the graphs as adjacency relations between parts}
	\label{Sierpinski-carpet}
\end{figure} 
\begin{figure}[h]
	\centering
	\includegraphics[width=0.9\textwidth]{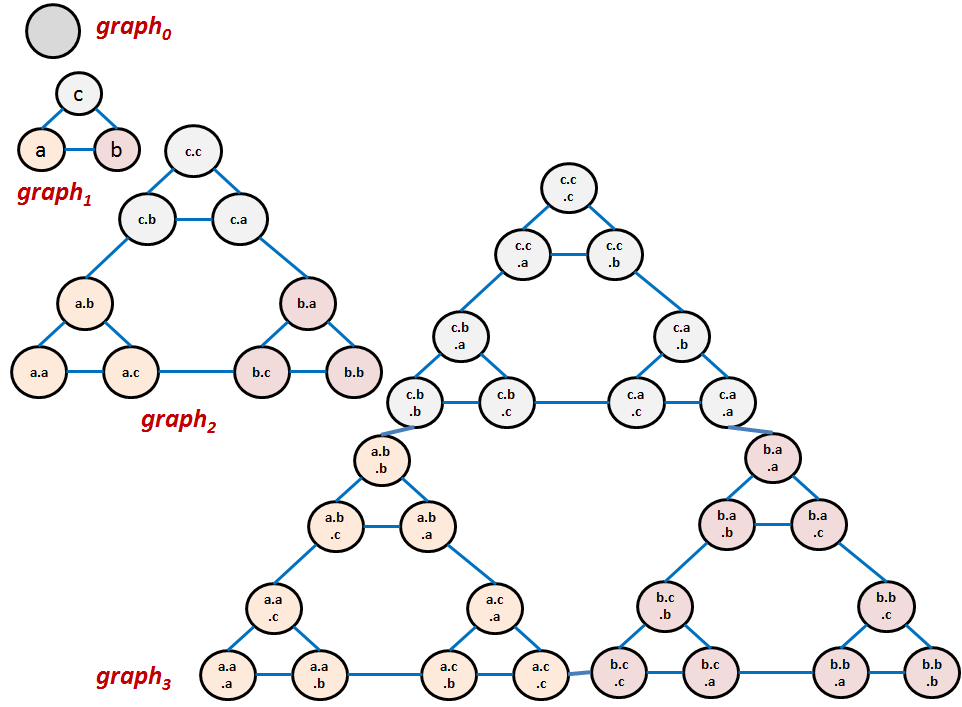}
	\caption{Consecutive partitions of an abstract Continuum corresponding to the Sierpinski sieve, and the graphs as adjacency relations between parts}
	\label{Sierpinski-sieve}
\end{figure} 
\begin{figure}[h]
	\centering
	\includegraphics[width=0.5\textwidth]{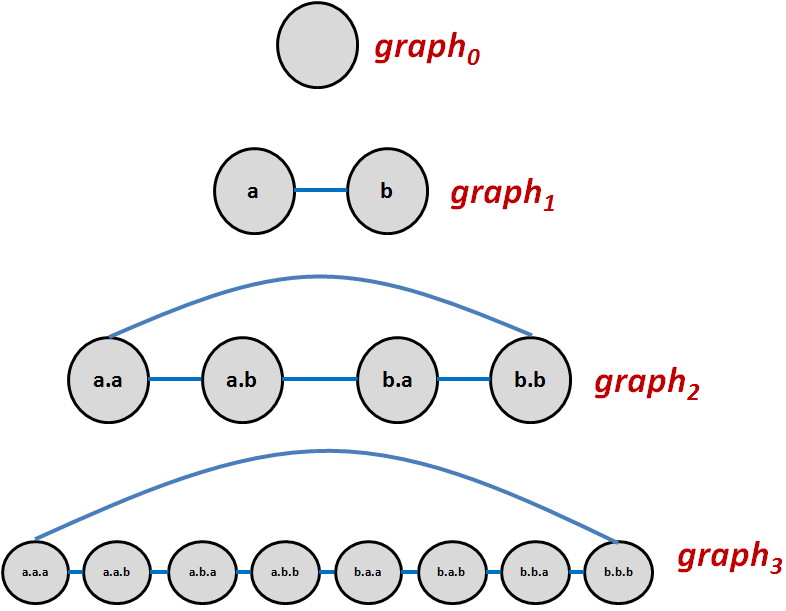}
	\caption{Circle. Edge (between the two border nodes) is added to each of the graphs }
	\label{circle}
\end{figure} 

\begin{figure}[h]
	\centering
	\includegraphics[width=0.91\textwidth]{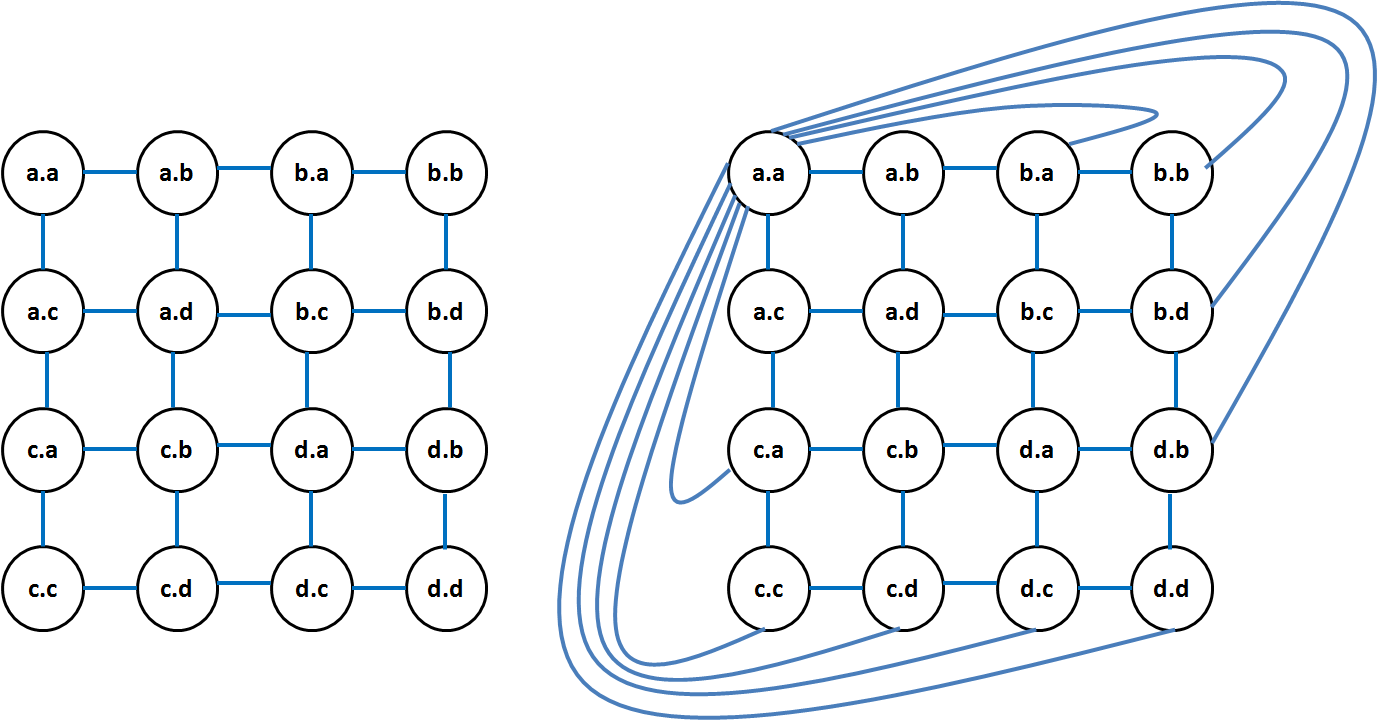}
	\caption{Graphs of the second division of the patterns corresponding to 2-cube, and 2-sphere (with one singular point) where edges are added between one fixed corner (e.g. nodes of the form $a.a.a.\ \dots \ a$) and the rest of the border nodes}
	\label{kwadrat-sfera}
\end{figure} 

\begin{figure}[h]
	\centering
	\includegraphics[width=0.99\textwidth]{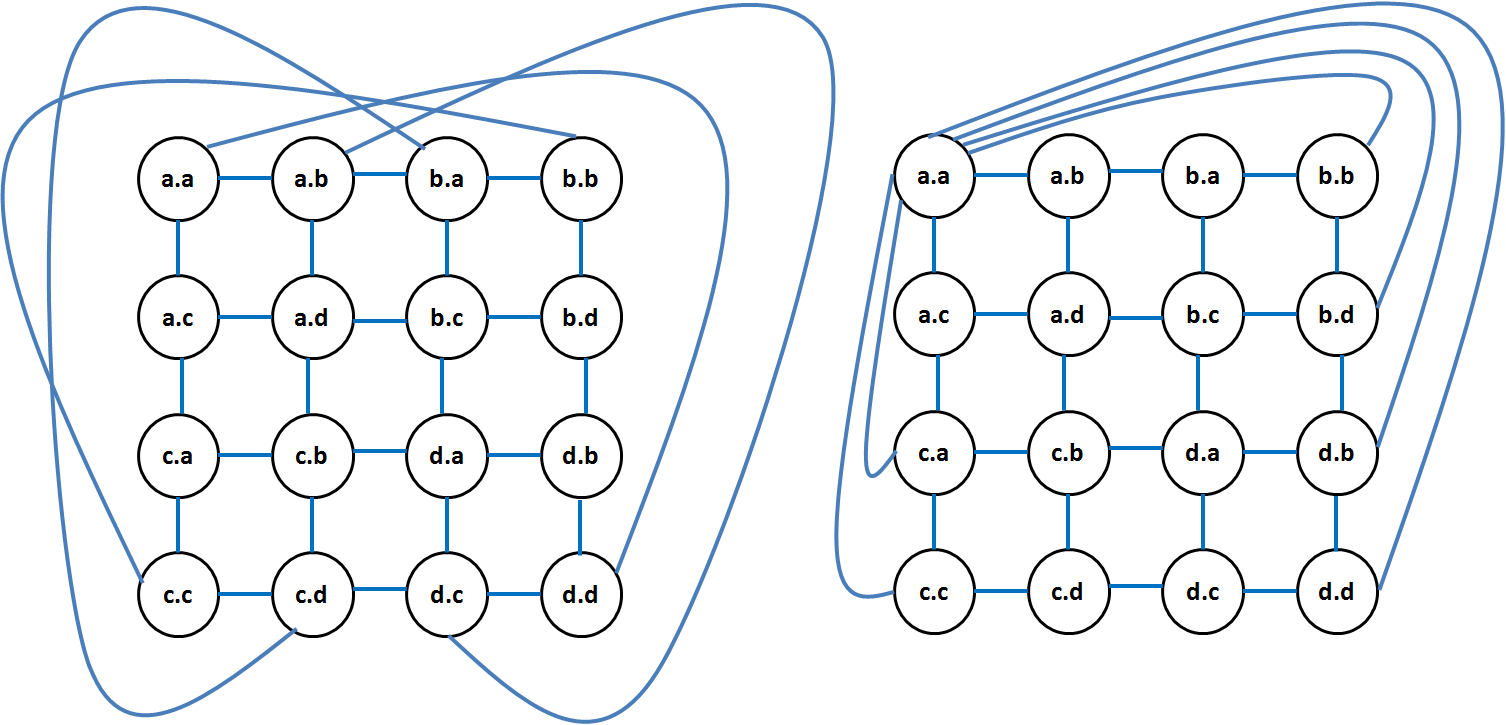}
	\caption{Graphs of the second division. The first from the left   corresponds to M\"{o}bius strip. If the two graphs are combined into one, then they correspond to the Klein bottle  with one singular point}
	\label{wstega-butelka}
\end{figure} 
\begin{figure}[h]
	\centering
	\includegraphics[width=0.65\textwidth]{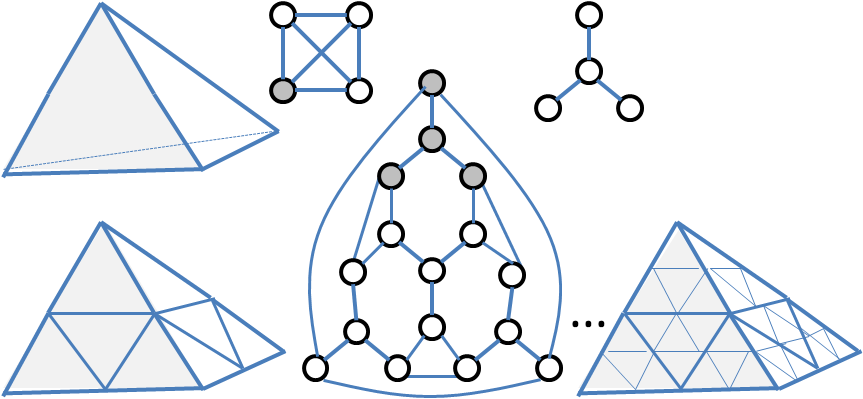}
	\caption{Tetrahedron surface. Note that any $graph_k$ has exactly four cycles of length 3 corresponding to the corners. The rest of the elementary cycles are of length 6 }
	\label{sfera-czworoscian}
\end{figure} 
\begin{figure}[h]
	\centering
	\includegraphics[width=0.6\textwidth]{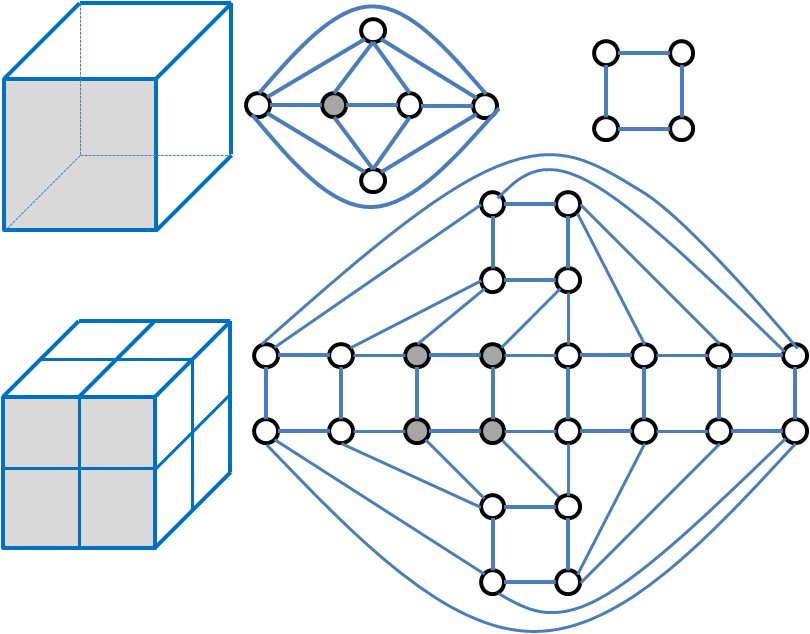}
	\caption{Hexahedron (3-cube) surface. Note that any $graph_k$ has exactly eight cycles of length 3 corresponding to the corners. The rest of the elementary cycles are of length 4 }
	\label{sfera-szescian}
\end{figure} 
In order to abstract completely from the intuitions of the unit cube, as well as from its interpretation as a subset of $\mathbb{R}^n$, only the above aspects should be taken into account in further investigations.  

The adjacency relations between parts in consecutive divisions may be represented as graphs. 
The examples (see Figures from \ref{unit-interval} to \ref{sfera-szescian}) show the divisions for Continua, and the corresponding adjacency relations. 

Let us take (for a while) the classic view and consider again the unit interval as Continuum. 
The interval may be represented as the $graph_0$, see Fig. \ref{unit-interval}, consisting only of a single node. 
The unit interval can be divided into two adjacent parts (two Continua). Each of the parts is the same (modulo size) as the original Continuum (the unit interval). This may be represented as the $graph_1$ consisting of two nodes, and one edge between them, corresponding to the adjacency of these two parts. 
We may continue the divisions, so that we get $graph_2$, $graph_3$, and so on, representing the parts as nodes, and adjacency as edges. 
Notice the  pattern of inductive construction of the sequence of graphs. 

Now, let the unit square be considered as the Continuum, see Fig. \ref{unit-square}. 
The method of construction of the sequence of graphs resulting from consecutive division of the square is similar. That is, the nodes represent small squares (Continua), whereas the edges represent the adjacency between them. 
The pattern of inductive construction of the sequence of graphs is important. 

Continuum may be different than Euclidean unit interval and unit square, see Fig. \ref{Sierpinski-carpet}. Let us consider the single node of the $graph_0$ that represents an abstract Continuum. 
The $graph_1$ is constructed from the $graph_0$, by dividing the single node (Continuum) into 8 nodes (Continua), and determining adjacency between the nodes. 
The $graph_2$ is constructed from the $graph_1$ by dividing each of the 8 nodes again into 8 nodes with the same adjacency (edges) as in the $graph_1$, and then some additional adjacency (edges) are determined. Each of the additional edges in $graph_2$ corresponds to an edge in the $graph_1$. 
Notice the pattern of inductive construction of the graphs. 

Similar inductive construction of graphs is presented in Fig. \ref{Sierpinski-sieve}. It corresponds to another fractal called the Sierpinski sieve.  
Notice the pattern of inductive division and adjacency. 

It seems that the sequence $(graph_0,\ graph_1,\ graph_2, \dots  graph_k, \dots  )$, and the pattern of its inductive construction, is crucial to grasp the proper grounding of the notion of Continuum. 

The patterns corresponding to the unit interval and the unit square are formally defined in Section \ref{Euclidean patterns}. The pattern for Sierpinski sieve is defined in Section \ref{fractals}. 

If more edges are added during the construction of the graph sequences for the unit interval and the unit square, then this results in some interesting Continua like circle in Fig. \ref{circle},  sphere in Fig. \ref{kwadrat-sfera} with one point singularity, and M\"{o}bius strip and Klein bottle (also with one singular point) in Fig. \ref{wstega-butelka}.  

Patterns corresponding to 2-sphere without singular points are presented in Fig. \ref{sfera-czworoscian} and Fig. \ref{sfera-szescian}. 
The patterns determine the same topological space,see Section \ref{ab}. 

The above examples do not exhaust inductive constructions of  sequences of the form  \\ 
$(graph_0,\  graph_1,\ graph_2,\ ...\ graph_k,\ ...\ )$ that determine interesting topological spaces like 3-cube, and in general, n-cubes, n-sphere, n-torus, and Riemannian manifolds in general. 

Let us summarize the above examples. 
Sequence $(graph_0,\  graph_1,\ graph_2,\ ...\ graph_k,\ ...\ )$ of the adjacency graphs constructed inductively according to a {\em division and adjacency pattern}, may be seen as the grounding for a Continuum. We are going to prove that.   


\section{ A  generic framework:  introducing a new primitive type corresponding to the Continuum  } 
\label{primitive-type}

What are types? There are primitive types, and there are types constructed from the primitive ones by using constructors like product, disjoin union, arrow, dependent type constructors $\Sigma$ and $\Pi$, and much more. 

Introducing a new primitive type requires: primitive objects, object constructors and destructors, and primitive relations.

Let $c_0$ denote an abstract unit as a primitive object. It corresponds to the $graph_0$ from the above examples.   
The object $c_0$ can be divided into parts. Each of the parts is again divided, and so on, inductively, to construct the sequence  $(graph_0,\  graph_1,\ graph_2,\dots graph_k,\dots )$. 

Each division is done according to the same division and adjacency pattern denoted by $(C_d, Adj_d
)$ where $C_d$ is a finite collection of elementary parts (denoted by $a_1,\ \ a_2,\dots, \ a_{l-1}, \ \ a_l$), \ \   $Adj_d$ is an adjacency relation between these parts
. The adjacency relation (as $graph_1$ from the examples) is connected, so that uniting all adjacent parts results again in the unit $c_0$. 

To denote the parts from consecutive divisions, it is convenient to use the {\em dot} notation introduced in Abstract Syntax Notation One (ASN.1). The notation has been already used in the above examples. So that, the parts resulted from division of $c_0$ are denoted by $c_0.a_1$,\ \  $c_0.a_2$, \dots , $c_0.a_{l-1}$, \ \ $c_0.a_l$. Each of the parts ($c_0.a_i$, for  $i=1, \ 2,  \dots  ,l$) is divided again, resulting in $c_0.a_i.a_1$,\  \  $c_0.a_i.a_2$,  \dots  , $c_0.a_i.a_{l-1}$, \  \ $c_0.a_i.a_l$. Analogously, it is done for the consecutive  divisions. 

The parts are objects of type (we are going to introduce) denoted by $C$. Actually, this is a simple example of Martin-L\"{o}f's W-type, that is, the type of lists (finite sequences) over $C_d$. However, contrary to W-types, for the type $C$ also a specific primitive adjacency relation will be introduced.

Since a formal introduction of primitive operations and primitive relations is not necessary to follow and understand the proposed grounding of the Continuum, it will be omitted here. It is done elsewhere, see \cite{TO}, where the importance for type theory is stressed and explained in detail. 
 
Let $C_k$ denote the collection of objects of length $k+1$ in $C$. There is a notational confusion with $C_d$ that denotes the  collection of elementary parts; the index $d$ is used only in this very context.  

\subsection{Adjacency as a primitive relation }
 \label{Adjacency} 

We are going to extend the adjacency relation  $Adj_d$ to all objects of the same length, i.e. to the objects in $C_k$ for arbitrary natural number $k$. This very extended relation corresponds to the $graph_k$ from the examples in  Section \ref{examples}. 

For a fixed $x$ of type $C$, adjacency between $x.a_1$,\ \ $x.a_2$,\ ..., \ $x. a_l$ (where $a_1$,\ $a_2$,\ ..., \ $a_l$ are the all elements belonging to $C_d$) is determined by $Adj_d$, that is, $x.a_i$ and  $x.a_j$ are defined as adjacent if $Adj_d(a_i, a_j)$ holds. So that, $(C_d, Adj_d)$ determines the division pattern, and it is called {\em d-pattern}. 

In order to extend the adjacency to all objects of $C_k$ (then, also to $C$), a merging pattern (in short {\em m-pattern}) is necessary to determine if there is adjacency between $x.a_i$ and $y.a_j$, where $x$ and $y$ are different adjacent objects of the same length. 

{\bf Important restriction on the m-pattern.} For any $x$ and $y$ of the same length, if they are adjacent, then there are $a_i$ and $a_j$ such that  $x.a_i$ and  $y.a_j$ are also adjacent. The above restriction follows from our intuition that if two parts are adjacent then after their divisions, there must be at least two sub-parts (one sub-part from each division) that are also adjacent.  


\subsection{Euclidean patterns} 
 \label{Euclidean patterns} 
 
For the unit interval $[0; 1]$, the d-pattern $(C_d, Adj_d)$ consists of two elements (parts), denoted by $a$ and $b$, that are adjacent. The consecutive three divisions are as follows, and they are also shown in Fig. \ref{unit-interval}: 
\begin{enumerate}
\item 
$c_0.a$ and $c_0.b$
\item 
$c_0.a.a$, \ \ $c_0.a.b$, \ \ $c_0.b.a$, \ \ $c_0.b.b$ 
\item 
$c_0.a.a.a$, \ \  $c_0.a.a.b$, \  \ $c_0.a.b.a$, \  \ $c_0.a.b.b$, \ 
 \  $c_0.b.a.a$, \  \ $c_0.b.a.b$,\ \ $c_0.b.b.a$, \ \ $c_0.b.b.b$ 
\end{enumerate}
The lexicographical order can be defined on the type $C$, if it is supposed that $a$ is preceding (lesser than) $b$.  

Since $c_0$ is the common prefix of all objects of type $C$, let it be omitted. 
For the first partition, $a$ and  $b$ are adjacent. 
For the second partition, also $a.b$ and  $b.a$ are adjacent.
According to the d-pattern, for any $x$ of type $C$, objects $x.a$ and  $x.b$ are adjacent. Also 
$x.a.b$ and  $x.b.a$ are adjacent.   

For the third partition, also  
$a.b.b$ and  $b.a.a$ are adjacent.

For the fourth partition, also 
$a.b.b.b$ and $b.a.a.a$, \ \ 
$a.a.b.b$ and $a.b.a.a$, \ \ 
$b.a.b.b$ and $b.b.a.a$  \ \ 
are adjacent. 

The m-pattern is is defined as follows. 
For any $x$ and $y$ of type $C$ of the same length, if $x.b$ and $y.a$ are adjacent, then  
\begin{itemize} 
\item if $x.b$ is less (in the lexicographical order) than $y.a$, then $x.b.b$ and $y.a.a$ are adjacent, 
\item  if $y.a$ is less than $x.b$, then $y.a.b$ and $x.b.a$ are adjacent. 
\end{itemize}
This completes the definition of the division and adjacency pattern for the unit interval. 

For the unit cube of dimension $2$ (i.e. the unit square), $C_d$ consists of four elements denoted by  $a, b, c, d$. Consecutive partitions are shown  in Fig. \ref{unit-square}.  The adjacency relation $Adj_d$ is between $a$ and $b$, $b$ and $d$, $d$ and $c$, and between $c$ and $a$. For any $x$ the following pairs of objects are adjacent: $x.a$ and $x.b$, $x.b$ and $x.d$, $x.d$ and $x.c$, $x.c$ and $x.a$. The m-pattern for dimension $2$ is based on our intuition of the 2-dimensional unit cube. 

More adjacent pairs of parts after the second and third division are shown in Tab. \ref{tab1} and Tab. \ref{tab2}.  
\begin{table}[h!]  
  \centering
   \begin{tabular}{c|c}
 $x.a.b$ and  $x.b.a$ & $x.a.d$ and $x.b.c$\\
 \hline
 $x.a.c$ and $x.c.a$ & $x.a.d$ and  $x.c.b$\\
 \hline
 $x.b.c$ and  $x.d.a$ & $x.b.d$ and  $x.d.b$\\
 \hline
 $x.c.b$ and $x.d.a$ & $x.c.d$ and $x.d.c$
  \end{tabular}
\caption{For any $x$ of $C$ the pairs in the table are adjacent}
\label{tab1}
\end{table}
\begin{table}[h!]
  \centering
   \begin{tabular}{c|c||c|c}
$x.a.b.b$ and $x.b.a.a$  &
$x.a.b.d$ and $x.b.a.c$  
&
$x.a.d.b$ and $x.b.c.a$ &
$x.a.d.d$ and $x.b.c.c$ \\
 \hline
$x.a.d.c$ and $x.c.b.a$ &
$x.a.d.d$ and $x.c.b.b$ 
&
$x.a.c.c$ and $x.c.a.a$ &
$x.a.c.d$ and $x.c.a.b$\\
 \hline
$x.b.c.c$ and $x.d.a.a$ &
$x.b.c.d$ and $x.d.a.b$ 
&
$x.b.d.c$ and $x.d.b.a$ &
$x.b.d.d$ and $x.d.b.b$ \\
 \hline
$x.c.b.b$ and $x.d.a.a$ &
$x.c.b.d$ and $x.d.a.c$ 
&
$x.c.d.b$ and $x.d.c.a$ &
$x.c.d.d$ and $x.d.c.c$ 
  \end{tabular}
\caption{Additional adjacent pairs after third division}
\label{tab2}
\end{table}

\begin{figure}[h]
	\centering
	\includegraphics[width=0.5\textwidth]{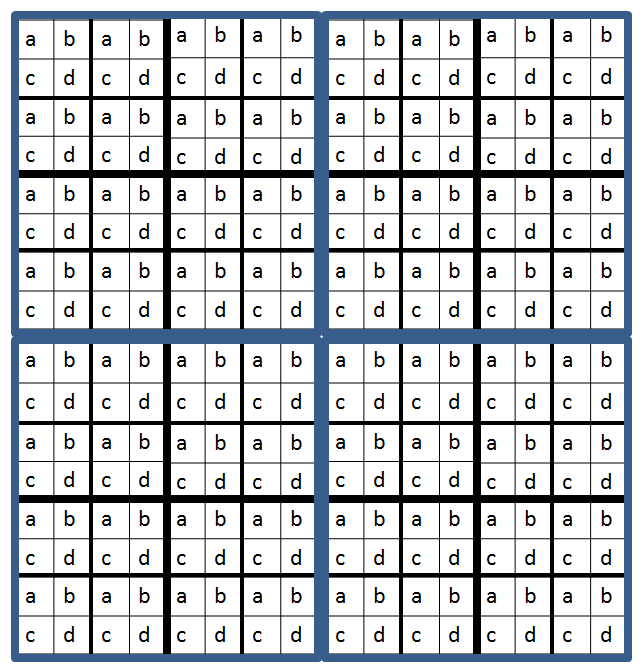}
	\caption{Unit square after 4-th division. The parts are denoted only by single letters $a$, $b$, $c$ and $d$.  Their complete denotation is of the form $o_1.o_2.o_3.o_4$, where $o_4$ is either $a$ or $b$ or $c$ or $d$, and $o_1$, $o_2$ and $o_3$  can be decoded (as one of the letters $a$, $b$, $c$ and $d$) from the position of the part in the unit square }
	\label{4square}
\end{figure} 
Let the alphabetical order on $a, b, c, d$ be extended to the lexicographical order on $C_k$ for any $k$.  
The m-pattern for the unit cube of dimension $2$ is defined inductively as follows. 
Let $x$ and $y$ be of the same length.  There are four following possible cases for which the adjacency relation is defined inductively. 
\begin{itemize}
\item  If $x.a$ and $y.b$ are adjacent and 
 \begin{itemize}
 \item $x.a$ is less than $y.b$ (according to the lexicographical order), then   
 $x.a.b$ and $y.b.a$ are adjacent, and 
  $x.a.d$ and $y.b.c$ are adjacent. 
 \item $y.b$ is less than $x.a$, then   
 $x.a.a$  and $y.b.b$ are adjacent, and 
 $x.a.c$ and $y.b.d$ are adjacent. 
 \end{itemize}
\item  If $x.b$ and $y.d$ are adjacent and 
 \begin{itemize}
 \item $x.b$ is less than $y.d$, then   
 $x.b.c$ and $y.d.a$ are adjacent, and 
  $x.b.d$ and $y.d.b$ are adjacent. 
 \item $y.d$ is less than $x.b$, then   
 $x.b.a$  and $y.d.c$ are adjacent, and 
 $x.b.b$ and $y.d.d$ are adjacent. 
 \end{itemize}  
\item  If $x.d$ and $y.c$ are adjacent and 
 \begin{itemize}
 \item $x.d$ is less than $y.c$, then   
 $x.d.b$ and $y.c.a$ are adjacent, and 
  $x.d.d$ and $y.c.c$ are adjacent. 
 \item $y.c$ is less than $x.d$, then   
 $x.d.a$  and $y.c.b$ are adjacent, and 
 $x.d.c$ and $y.c.d$ are adjacent. 
 \end{itemize}
\item  If $x.a$ and $y.c$ are adjacent and 
 \begin{itemize}
 \item $x.a$ is less than $y.c$, then   
 $x.a.c$ and $y.c.a$ are adjacent, and 
  $x.a.d$ and $y.c.b$ are adjacent. 
 \item $y.c$ is less than $x.a$, then   
 $x.a.a$  and $y.c.c$ are adjacent, and 
 $x.a.b$ and $y.c.d$ are adjacent. 
 \end{itemize}    
\end{itemize}  
The Fig. \ref{4square} may help to verify the above cases. 
This completes the definition of the Euclidean pattern for the unit square. 

Note that the above m-pattern is not so simple, so that in order to grasp it we must refer to the pictures, i.e. to our intuitions. 

In the very similar way, for any $n$ dimensional cube, the appropriate d-pattern and m-pattern can be formally defined. Let it be called {\em the Euclidean pattern of dimension $n$}. 

\section{Generic form of the patterns} 
\label{generic form}

In its generic form, the division and adjacency pattern inductively determines the adjacency relation (graph) on the objects (nodes) of $C_k$ for any $k$ that corresponds to $graph_k$ from the examples in the Section \ref{examples}. 
From now on, d-pattern may be arbitrary, and is determined by a collection $C_d$ of elementary parts (primitive objects) $a_1, \ a_2, \ \dots \ a_l$; and an arbitrary relation $Adj_d$ (a connected graph) on these parts.  

 m-pattern  that (along with d-pattern) determines the adjacency relation (denoted by $Adj_{C_k}$) for objects of $C_k$, for arbitrary $k$, also may be arbitrary. The only necessary condition is that the graph related to $Adj_{C_k}$ must be connected. Note that $Adj_{C_k}$ is symmetric, i.e. for any $x$ and $y$ in $C_k$, if $Adj_{C_k}(x;y)$, then $Adj_{C_k}(y;x)$. It is also convenient to assume that $Adj_{C_k}$ is reflexive, i.e. any $x$ is adjacent to itself, formally $Adj_{C_k}(x;x)$. 
 
The relation $Adj_{C_k}$ can be extended (by appropriate construction) to the type  $C$, that is, to the relation $Adj_C$ in the following way. For any $x$ and $y$ of different length (say $n$ and $k$ respectively, and  $n<k$), relation $Adj_C(x;y)$ (as well as its symmetric version $Adj_C(y;x)$) holds if there is object $\bar{x}$ of length $k$ such that $x$ is a prefix (initial segment) of $\bar{x}$,  and $Adj_{C_k}(\bar{x}; y)$ holds, i.e.  $\bar{x}$ is of the same length as $y$ and is adjacent to $y$.  


Let us note that the above definitions of d-pattern and m-pattern are so general that the resulting adjacency  relation, i.e. the sequence $(graph_0, graph_1, graph_2 , \dots graph_k,  \dots  )$, may be  different than the adjacency relation for the unit cubes. 
The only necessary condition imposed on the patterns is that the corresponding adjacency relation on $C_k$, as $graph_k$, must be connected for any $k$. That is, the result of uniting all adjacent parts must be again the initial unit object $c_0$.  

 
Let d-pattern and m-pattern together be denoted by $\mathbb{P}$. 
For a fixed pattern $\mathbb{P}$, let $C$ be denoted by $C_{\mathbb{P}}$ whereas $Adj_C$ by $Adj_{\mathbb{P}}$. 
Analogously,
$C_{\mathbb{Q}}$ denotes the primitive type for pattern $\mathbb{Q}$, and $Adj_{\mathbb{Q}}$ denotes its adjacency relation.  Euclidean pattern of dimension $n$ is denoted by $\mathbb{E}^n$.  

Let us recall the intuitive meaning of Continuum: {\em `` ... no part of which can be distinguished from neighboring parts except by arbitrary division''}. That is, the self-similar fractal structure resulting from the adjacency patterns is not enough. The following conditions are essential.   
\begin{itemize}
\item  {\em  Indiscernibility of the parts:}  
The parts from $C_d$ can not be distinguished only on the basis of their neighborhood (adjacency) structures. That is, any part in $C_d$ has exactly the same neighborhood structure. 
%
\item {\em Homogeneity of the parts:}  For any $x$ and $y$ (except the border objects) in $C_k$, their neighborhoods in $C_k$ are isomorphic with regard to the adjacency relation, i.e.  any object from  $C_k$ (except the border objects) has exactly the same fixed adjacency structure in $C_k$.   
%
\end{itemize}
Note that Sierpinski carpet (see Fig. \ref{Sierpinski-carpet}) satisfies only the Indiscernibility condition. 
Sierpinski sieve (see Fig. \ref{Sierpinski-sieve}) and the Euclidean patterns satisfy the Indiscernibility condition and Homogeneity condition. 
The natural question is if there are other patterns that satisfy these two conditions.

\subsection{Patterns and fractals} 
\label{fractals}

Examples of non-Euclidean patters are shown (in the form of graph sequences) in Fig. \ref{Sierpinski-carpet}, Fig. \ref{Sierpinski-sieve}, and Fig. \ref{circle}. For the Sierpinski sieve, the d-pattern consists of three parts (denoted by $a$, $b$, and $c$), and any two of these parts are adjacent. The inductive definition of m-pattern is as follows. 
For any $x$, the objects $x.a$, $x.b$, and $x.b$ are pairwise adjacent; and also $x.a.b$ and $x.c.b$ are adjacent, $x.b.a$ and $x.c.a$ are adjacent, $x.a.c$ and $x.b.c$ are adjacent. 
\\ 
Also $x.a.y$ and $x.c.y$ are adjacent, if $y$ is a suffix consisting of elements that all are $b$ (that is, $b.b.\  \dots b$). 
\\ 
Objects $x.b.v$ and $x.c.v$ are adjacent, if $v$ is a suffix consisting of elements that all are $a$. 
\\ 
Objects $x.a.z$ and $x.b.z$ are adjacent, if $z$ is a suffix consisting of elements that all are $c$. 

The above example corresponds to the famous fractal called Sierpinski triangle, see \cite{Sierpinski}, and  Tower of Hanoi graph, see \cite{hinz2017survey}. However there, the adjacency relation is defined as {\em ``adjacent-endpoint relation''}, see \cite{lipscomb2008fractals}. 
That is, an infinite sequence $a_1.a_2, \dots$ (of elements form finite set $A$)  is called an endpoint
if there exists index $k$ such that $a_{k+n} = a_k$ for all $n$.  Two distinct
endpoints $a_1.a_2, ...$ and $b_1.b_2, ...$ are adjacent if there is $k$ such that for all $i$ less than $k$: $a_i = b_j$ (i.e. the endpoints prefixes of length $k-1$ are the same), and there are different $a$ and $b$ in $A$ such that $a_k$ is $a$ and $b_k$ is $b$, and for all $j$ greater than $k$: $a_j$ is $b$ and $b_j$ is $a$. 
The adjacent-endpoint relation induces a fractal structure. 

The primitive types of the general form $C_{\mathbb{P}}$ may be seen as a generalization of this classical definition of fractal structure. 

Another abstract definition of fractal structures refers to iterated function systems, see \cite{falconer2004fractal}. A system is a finite collection of contraction mappings on a complete metric space $(X, \rho)$, say $f_i: X \rightarrow X$  (for $i \leq k$) such that for any $i$ there is a positive real number $r_i < 1$ such that for any $x_1$ and $x_2$ in $X$: $\rho(f_i(x_1), f_i(x_2)) \ < \ r_i \rho(x_1, x_2)$. Then, the unique subset $F$ (called attractor or fractal) of $X$ is defined by the fix-point equation $F = \bigcup^k_{i=1}f_i(F)$. The fractal results as the limit of iterative applications of the functions $f_i$, i.e. $F= \dots \bigcup^k_{i=1} f_i((\bigcup^k_{i=1} f_i(\bigcup^k_{i=1}f_i(X))))$.  

The constructions introduced in the book \cite{lipscomb2008fractals} 
are of particular interest here. They are based of the adjacent-endpoint relation of topological spaces corresponding to 2-sphere, and are very similar to the one presented in this paper. However, still the adjacency relations on finite sequences (graph nodes) does not appear explicitly in the work of Lipscomb. There is also interesting correspondence (presented in this book) between iterated  function systems and the adjacent-endpoint relation.  

Finite version of adjacent-endpoint relation for Sierpinski graphs was considered by \cite{Milutinovic}. Actually, it leads to a specific sequence of graphs $(graph_0, graph_1, graph_2 , \dots graph_k,  \dots  )$, so that it can be identified with a specific division and merging pattern. 
Extension of adjacent-endpoint relation to generalized Sierpinski graphs was proposed in  \cite{hinz2017survey} and \cite{klavzar2017connectivity}. 
The difference is that for Sierpinski graphs, the basis is the complete graph $K_l$ for a fixed number $l$ of nodes, whereas for the generalized Sierpinski graph arbitrary connected graph may be taken as the basis.  

Sierpi{\'n}ski Carpet Graph, see Fig. \ref{Sierpinski-carpet}, is also  interesting. Its common definition (in a sense, also a construction) invokes the unit square $[0; 1]^2$ a subset of $\mathbb{R}^2$, see \cite{cristea2010connected}. The inductive construction consists in the nested procedure of dividing square into 9 smaller squares and removing the one in the middle. Similar construction may be applied to obtain the Sierpinski sieve. Both fractals are generalization of the construction of Cantor set to the two dimensional Euclidean space. Each of the fractals (as a subset of Euclidean space) has a natural topology induced by the topology of that Euclidean space.

Note that W-types augmented with the adjacency relations (determined by division and adjacency patterns $\mathbb{P}$) may serve as a universal method for fractal construction without using abstract notions like complete metric spaces. 
Actually, topologies generated by the patterns $\mathbb{P}$ may be   different and independent of the Euclidean topologies, see the next Section \ref{ab}.  

The d-pattern and the m-pattern corresponding to the Sierpi{\'n}ski Carpet Graph are easy to construct.

The above examples clearly show that Continua as types, with adjacency relations constructed by non-Euclidean patterns, are also interesting. 


Let us summarize this section. For any d-pattern and any m-pattern (that together determine division and adjacency pattern $\mathbb{P}$), a primitive W-type can be introduced (denoted by $C_{\mathbb{P}}$) along with adjacency relation $Adj_{\mathbb{P}}$. Usually, the type corresponds to a fractal.  However, pattern $\mathbb{P}$ of inductive construction of sequence $(graph^{\mathbb{P}}_0, graph^{\mathbb{P}}_1, graph^{\mathbb{P}}_2 , \dots  graph^{\mathbb{P}}_k,  \dots )$ may be arbitrary, and sometimes it does not have a fractal-like structure. 

Let us note that a sequence of the form $(graph_0, graph_1, graph_2 ,  \dots  graph_k,  \dots  )$ and corresponding adjacency relation may be constructed inductively not necessary using  d-pattern and  m-pattern. Hence,  pattern $\mathbb{P}$ is defined as an inductive method for generating the sequence  $(graph^{\mathbb{P}}_k, \ k= 0,1,2,  \dots  )$ and corresponding adjacency relation $Adj_{\mathbb{P}}$. 

\section{Abstraction from patterns to topological spaces } 
\label{ab}

The primitive types augmented with adjacency relations correspond to important and interesting topological spaces as was shown in Section \ref{examples}. It was claimed there that such types are the grounding for the spaces that themselves are abstract notions. Now, we are going to show how to abstract from a primitive type to a topological space. 
 
Functions
on the types (say from $C_{\mathbb{P}}$ into $C_{\mathbb{Q}}$),  are important for the abstraction. The patterns $\mathbb{P}$ and $\mathbb{Q}$ are arbitrary. 

A function $f:C_{\mathbb{P}} \rightarrow C_{\mathbb{Q}}$ is defined as {\em  monotonic} if for any  $x$ and $y$ such that $x$ is  prefix of $y$,  $f(x)$ is a prefix of $f(y)$, or $f(x)$ is equal to $f(y)$. 

Let us consider sequences of the form $(x_k;\  k=1,2, \dots)$ of objects of type $C_{\mathbb{P}}$ such that $x_k$ is of length $k$ and $x_k$ is a prefix of $x_{k+1}$.  Monotonic function $f$ is {\em strictly} monotonic if for any such sequence $(x_k;\  k=1,2, \dots)$ and for any $n$ there is  $j > n$ such that $f(x_n)$ is not $f(x_j)$.

A strictly monotonic function $f: C_{\mathbb{P}} \rightarrow C_{\mathbb{Q}}$ is defined as {\em continuous} if for any two adjacent objects $x$ and $y$, the objects $f(y)$ and $f(x)$ are also adjacent. 
This continuity may be seen as more intuitive if $x$ and $y$ are of the same length. Actually, it is the same; it follows from the definition of adjacency relation for objects of different length (see Section  \ref{generic form}), and from the strict monotonicity.  

The strictly monotonic and continuous functions are important for the Brouwer's Continuum and the famous Brouwer's Continuity Theorem, as we will see in the next sections. 

Let $C_d$ correspond to pattern $\mathbb{P}$, and consider the type $C_{\mathbb{P}}$. 
Let us introduce the abstract notion of infinite sequences over $C_d$ and denote the set of such infinite sequences by  
$C^\infty_{\mathbb{P}}$. For an infinite sequence denoted by $u$, let $u(k)$ denote its prefix of length $k$, i.e. the initial finite sequences of length $k$. 

The set $C^\infty_{\mathbb{P}}$ may be considered as a topological space (Cantor space) with topology determined by the family of open sets $U_x$ such that for any $x$ of length $n$:  $U_x = \{u: u(n) = x\}$. Note that the pattern $\mathbb{P}$ is not used in the definition.  

We are going to introduce a different topological structure on  the set $C^\infty_{\mathbb{P}}$ determined by the pattern. In order to do so the pattern $\mathbb{P}$ must satisfy the following two properties. 
\\
{\bf Property 1.} {\em If $x$ and $y$ (of the same length) are adjacent, then there are two elements of $C_d$ (say $c_1$ and $c_2$) such that $x.c_1$ and $y.c_2$ are adjacent}. 
\\ 
{\bf Property 2.}  {\em If $x$ and $y$ (of the same length) are not adjacent, then for any two elements $c_1$ and  $c_2$ of $C_d$,  the objects $x.c_1$ and $y.c_2$ are not adjacent}. 
\\ 
Let a pattern that satisfies the above two conditions be called {\em regular}. For regular patterns, if $x$ and $y$ (of the same length $k$) are adjacent, then their prefixes of length $i<k$, i.e. $x(i)$ and $y(i)$, also are adjacent; they may be also the same. 

For a regular pattern $\mathbb{P}$, the corresponding $(graph^{\mathbb{P}}_k, \ k= 0,1,2,  \dots  )$ is an inverse sequence and  $C^\infty_{\mathbb{P}}$ is the well known inverse limit of the sequence, see   \cite{Alexandroff}, and \cite{debski2018cell} for a recent review of the subject. 

Euclidean patterns and their extensions (see Fig. \ref{circle}, Fig. \ref{kwadrat-sfera} and Fig.  \ref{wstega-butelka}) are regular. Also the patterns that correspond to Sierpinski sieve (see Fig. \ref{Sierpinski-sieve}) and Sierpinski carpet (see Fig. \ref{Sierpinski-carpet}) are regular.  

From now on, we consider only regular patterns. So that the term  {\em pattern} means regular pattern.  The reason for restricting patterns to the regular ones, is the following definition of extension of adjacency relation to the infinite sequences.  

Two infinite sequences $u$ and $v$ are defined as {\em adjacent} if for any $k$, the prefixes $u(k)$ and $v(k)$ are adjacent, i.e.  $Adj_{\mathbb{P}}(u(k) ; v(k) )$ holds. Let this adjacency relation, defined on the infinite sequences, be denoted by $Adj^\infty_{\mathbb{P}}$. Note that any infinite sequence is adjacent to itself. Two different adjacent infinite sequences may have the same prefixes.  

The transitive closure of $Adj^\infty_{\mathbb{P}}$ is an equivalence relation denoted by $\sim$. Let the quotient set be denoted by $C^\infty_{\mathbb{P}}/_{\sim}$, whereas its elements, i.e. the equivalence classes be denoted by $[v]$ and $[u]$.  The topological space $C^\infty_{\mathbb{P}}/_{\sim}$ is the  Hausdorff reflection of $C^\infty_{\mathbb{P}}$ relatively ti the equivalence relation ${\sim}$. 

Let us take the pattern $\mathbb{E}^1$, corresponding to the unit interval $[0; \ 1]$, as an example.  
Let us recall that $C_d$ consists of $a$ and $b$ that are adjacent.  $C_{\mathbb{E}^1}$ is the type of finite sequences over $C_d$. Objects of this type are denoted, in the dot notation, as $c_0.c_1.c_2.\ \dots  c_k$ where $c_i$ (for $i=1,2,  \dots ,k$) are elements of $C_d$. Since $c_0$ is the common prefix of all objects, it is omitted. 
For any $x$ of length $k$, let the infinite sequences $v$ and $u$ be such that $x\ =\ v(k)\ =\ u(k)$, and the $(k+1)$-th element of $v$ is $b$, and the $(k+1)$-th element of $u$ is $a$, and for any $n$ the $(k+1+n)$-th element of $v$ is $a$, and the $(k+1+n)$-th element of $u$ is $b$.  
That is, $v=x.b.a.a.a.\ \dots $ and $u=x.a.b.b.b.\ \dots $.  Then, $v$ and $u$ are adjacent, and the equivalence class $[v]$ consists exactly of two elements $v$ and $u$. For the rest of infinite sequences $z$, the equivalence classes $[z]$ are singletons, i.e. each one consists of one element only. 
 
The transitive closure of the adjacency relation $Adj^\infty_{\mathbb{E}^1}$ does not extend the relation; it remains the same. However, it is not so for higher dimensional patterns, e.g. $\mathbb{E}^2$, where some equivalence classes have 4 elements, others have 2 elements, and the rest classes are singletons. 

Let $\mathbb{P}$ be again an arbitrary regular pattern. The relation $Adj_{\mathbb{P}}$ determines natural topology on the quotient set $C^\infty_{\mathbb{P}}/_{\sim}$. 

Usually, a topology on a set is defined by the family of open subsets that is closed under finite intersections and arbitrary unions. The set and the empty set belong to that family.  Equivalently, the topology is determined by families of closed neighborhoods of the points of the set. For any equivalence class $[v]$ (i.e. a point in $C^\infty_{\mathbb{P}}/_{\sim}$), the neighborhoods (indexed by $k$) are defined as the sets $U^{[v]}_k$ of equivalence classes $[u]$ such that there are $v'\in [v]$ and $u'\in [u]$ such that $u'(k)$ is adjacent to $v'(k)$, i.e. $Adj_{\mathbb{P}}(u'(k);\ v'(k))$ holds. Note that $u'(k)$ may be equal to $v'(k)$. 

A sequence $(e_1, \ e_2,\ \dots \ e_n, \ \dots \ )$ of  equivalence classes of $\sim$ (elements of the set $C^\infty_{\mathbb{P}}/_{\sim}$) converges to $e$, if for any $i$ there is $j$ such that for all $k>j$, \ \ \  $e_k \in U^e_i$. 

Note that the topological space $C^\infty_{\mathbb{P}}/_{\sim}$ is is a Hausdorff  compact space. It is an abstraction done from the type $C_{\mathbb{P}}$ and its adjacency relation $Adj_{\mathbb{P}}$. The abstraction uses actual infinites, i.e. infinite sequences.   

For the pattern $\mathbb{E}^1$ corresponding to the unit interval $[0;\ 1]$, the space $C^\infty_{\mathbb{E}^1}/_{\sim}$ is homeomorphic to that interval. If the natural linear lexicographical order is introduced, then the space is isomorphic to the interval.  

For another patterns $\mathbb{P}$ from the examples in Section \ref{examples}, the space is homeomorphic to one of the topological spaces: the unit square, Sierpinski sieve, 2-sphere, etc.  Hence, the patterns and the resulting types should be considered as the groundings for these topological spaces. 

The patterns may be quite sophisticated so that the  abstracted topological spaces may be homeomorphic to important spaces like $n$-spheres, M\"{o}bius strip and Klein bottles of higher dimensions,  Sierpinski sieve and carpet, and many others fractals and Riemannian manifolds.   

The extension of a strictly monotonic function $f: C_{\mathbb{P}} \rightarrow C_{\mathbb{Q}}$ to the function $f^{\infty}$ from $C^{\infty}_{\mathbb{P}}$ into $C^{\infty}_{\mathbb{Q}}$ is natural. That is, $f^{\infty}(v)$ is defined as element $z$ in $C^{\infty}_{\mathbb{Q}}$ such that for any $i$ there is $k$ such that $z(i)$ is a prefix of $f(v(k))$. 
The strict monotonicity of $f$ is essential here.  An equivalent definition of the extension is that  $f^{\infty}(v)$ is the limit of the sequence ($f(v(n))$, for $n= 1,2, \ \dots $). By the monotonicity, for any $n$, $f(v(n))$ is a prefix of $f(v(n+1))$ or equal to $f(v(n+1))$. By the strict monotonicity, the length of consecutive elements of the sequence is unbounded. 

Let us consider again the examples from the Fig. \ref{sfera-czworoscian} (where the pattern is denoted by $\mathbb{P}$), and Fig. \ref{sfera-szescian} (pattern $\mathbb{Q}$). The both patterns correspond to the same topological space, i.e. to 2-sphere. 
There is no strictly monotonic and continuous function between $C_{\mathbb{P}}$ and $C_{\mathbb{Q}}$. 

\section{Brouwer's choice sequences}
\label{cs} 

Let us consider again the type $C_{\mathbb{E}^1}$ where the pattern $\mathbb{E}^1$ corresponds to the unit interval $[0;\ 1]$, that is, a subset of the real numbers. 

The finite sequences may be interpreted (according to Fig.  \ref{unit-interval}) as the subintervals of the consecutive partitions of the unit interval $[0;\ 1]$. That is, $a$ is $[0; \frac{1}{2}]$ and $b$ is  $[\frac{1}{2}; 1]$, \ \ $a.a$ is $[0; \frac{1}{4}]$,\ \  $a.b$ is $[\frac{1}{4}; \frac{2}{4}]$, \ \  $b.a$ is $[\frac{2}{4}; \frac{3}{4}]$, \ \ $b.b$ is $[\frac{3}{4}; 1]$, and so on. Then, the adjacency (according to pattern $\mathbb{E}^1$) between the sequences of length $k$ is determined by the adjacency of the sub-intervals resulting from the $k$-th division. However, it is a bit misleading because it is based on our intuition, whereas the aim is to provide the concrete grounding for this intuition. 
The pattern of this adjacency (without referring to this intuition) is simple, and is the basis for inductive construction of the adjacency relation for the objects of length $k$, and then to extend the relation  to $Adj_{\mathbb{E}^1}$ on all objects of type $C_{\mathbb{E}^1}$. This has been already done in the Section \ref{Euclidean patterns}. 

Now let us consider type $C_{\mathbb{P}}$ where the regular pattern $\mathbb{P}$ is arbitrary. 
For any infinite sequence $v$, the sequence of finite sequences ($v(n)$, for $n= 1, 2, \dots $) is defined as {\em a choice sequence} over $C_d$. Although from the formal point of view, a choice sequence is the same as the corresponding infinite sequence, according to Brouwer, it makes the difference, see the Intuitionism entry at Stanford Encyclopedia of Philosophy \url{https://plato.stanford.edu/entries/intuitionism/}. A choice sequence is a potential infinity, whereas an infinite sequence is an actual infinity. It means that in mathematical constructions, if we take a choice sequence, we use only its prefixes (finite sequences), whereas if we take infinite sequence, we use it as an actual infinity, as it was already done in the definition of topological space $C^\infty_{\mathbb{P}}/_{\sim}$. 

It is convenient to use the same symbols (e.g. $u$, $v$, $z$) to denote infinite sequences and choice sequences.  It will be clear, from the context, what they denote. 

Note that for the pattern $\mathbb{E}^1$ corresponding to the unit interval, the topological space $C^\infty_{\mathbb{E}^1}/_{\sim}$ is homeomorphic to the unit interval $[0; 1]$ with the usual Euclidean topology. 
Then, a choice sequence $v(n)$, for $n= 1, 2, \dots $ may be interpreted as a sequence of nested intervals converging to a point in $[0;\ 1]$. 

Note that the sequence of adjacency graphs $(graph^{\mathbb{P}}_k$, for $k=0,1,2,  \dots  )$ 
corresponds to  Brouwer's notion of spread, whereas the set of nodes of $graph^{\mathbb{P}}_k$ to a Brouwer's fun. However, the adjacency determined by $graph^{\mathbb{P}}_k$ is more rich structure than a tree with branches of length $n$. 
The notion of choice sequences as well as spread and fun can be defined for any pattern $\mathbb{P}$. However, the adjacency relation is not used in these definitions.  

The Axiom of Continuity is crucial for the Brouwer's approach to Continuum.
The weak version of the axiom is as follows. 
\\ {\em 
For any predicate $A(v, n)$ on the choice sequences $v$ and the natural numbers $n$:\\
if for any choice sequence $v$ there is $n$ such that $A(v, n)$, \\
then for any choice sequence $v$ there are natural numbers $m$  
and $n$ such that for any choice sequence $u$: \\ 
if  $u(m) = v(m)$, then $A(u, n)$.}

It means that, although the choice sequences are potentially infinite, the truth of the predicate $A(v, n)$ is determined by an initial segment (prefix) of $v$. The strong axiom of continuity requires that this very initial segment can be effectively determined by a computable functional taking a choice sequence as an argument, and returning the natural number as the length of the desired prefix of the argument.
 
In our framework, the weak axiom may be rewritten as follows.  
\\  
{\bf The postulate.} {\em Any function $F$ on choice sequences (i.e. from $C^\infty_{\mathbb{P}}$ into $C^\infty_{\mathbb{Q}}$) is the extension of a strictly monotonic function $f$ of type $C_{\mathbb{P}} \rightarrow C_{\mathbb{Q}}$}.  \\ 

The extension (defined at the end of Section \ref{ab}) implies that if $F(v) = f^\infty(v) = z$, then for any $i$ there is $k$ such that $z(i)$ is a prefix of $f(v(k))$. So that in order to compute an initial segment (prefix) of the value $F(v)$, the function $F$ needs only an initial segment of its argument $v$. So that, the computations done by the function $F$ are only potentially infinite. 

In the classic Analysis, there are functions that contradict the postulate. As an example, let us define function $F$ on $C^\infty_{\mathbb{E}^1}$ \ (homeomorphic to the unit interval $[0; \ 1]$) in the following way.  If $v$ corresponds to a rational number, then let $F(v)$ be the infinite constant sequence where each its element is $a$.  Otherwise, let the infinite sequence $F(v)$ be also constant and consists of elements that all are equal to $b$.  Another simple example is the function $G$ such that for any $v$ that corresponds to a real number less or equal $\frac{1}{2}$,  the infinite sequence $G(v)$ is constant and consists of elements that are equal to $a$.   Otherwise, $G(v)$ consists of elements that each one is $b$.   For the choice sequence $v$ such that $v(1) = b$, and $k$-th element of $v$ is $a$ for $k>1$, the determination of the value $G(v)$ must take the whole infinite sequence $v$. Note that the functions $F$ and $G$ are not continuous. 

The Postulate may be reformulated using predicates instead of functions in the following way. 
\\
Let $f$ be the strictly monotonic function corresponding to $F$ from the postulate.  The predicate $A_F$ is defined as follows. \\
$A_F(v,z,i)\ \equiv \ \exists_k \ z(i)$ is a prefix of $f(v(k))$. 
\\
The truth of the predicate $A_F(v,z,i)$ is determined by an initial segment (prefix) of $v$ and a prefix of $z$. 	Hence, the predicate satisfies The Axiom of Continuity.  
The predicate $\forall_i A_F(v,z,i)$ is true if and only if $f^{\infty}(v) \ = \ F(v)\ =\ z$. 

The Postulate is sufficient to prove the famous Brouwer Continuity Theorem. It is done in the Section \ref{btc}.

\section{Functions on types and their extensions to functions on topological spaces }

Let us consider again a strictly monotonic function $f: C_{\mathbb{P}} \rightarrow C_{\mathbb{Q}}$, and its extension to $f^\infty: C^\infty_{\mathbb{P}} \rightarrow C^\infty_{\mathbb{Q}}$. 

The transformation of $f^{\infty}$ to the function $f^{\infty}/_{\sim}$ from $C^{\infty}_{\mathbb{P}}/_{\sim}$ into $C^{\infty}_{\mathbb{Q}} /_{\sim}$ is not always possible. 

If $f$ is not continuous, then $f^{\infty}/_{\sim}$ need not to be a function. That is, there may be two adjacent infinite sequences $v$ and $u$ (i.e. $v \sim u$) such that $f^{\infty}(v)$ is different than $f^{\infty}(u)$. Since $[v]$ is the same as $[u]$, the value$f^{\infty}/_{\sim}([v])$ can not be defined.

For the type $C_{\mathbb{P}}$ and type $C_{\mathbb{Q}}$ where the both patterns $\mathbb{P}$ and $\mathbb{Q}$ are $\mathbb{E}^1$ (i.e Euclidean pattern of dimension 1 corresponding to the unit interval), the proof is simple. Suppose that $f$ is strictly monotonic and not continuous. Then, there are sequences $x$ and $y$ (of the same length, because $\mathbb{E}^1$ is a regular pattern) that are adjacent, and $f(x)$ and $f(y)$ are not adjacent. So that, either $x.a$ and $y.b$ are adjacent, or $x.b$ and $y.a$ are adjacent. That is, either $x.a.a.a.a.a.\ \dots $ and $y.b.b.b.b.b.\ \dots $ are equivalent (adjacent), or $y.a.a.a.a.a.\ \dots $ and $x.b.b.b.b.b.\ \dots $  are equivalent. Let $r_1$ and $r_2$ be the infinite sequences from the above ones that are equivalent, i.e. they represent the same real number. By the assumption that $f(x)$ and $f(y)$ are not adjacent, and by the strict monotonicity of $f$and the regularity of pattern $\mathbb{E}^1$, \ \  $f^{\infty}(r_1)$ and $f^{\infty}(r_2)$) are not equivalent and correspond to different real numbers. 
\\
Note that in the proof presented above, only the regularity of the pattern  $\mathbb{E}^1$ was essential. That is, the construction of $r_1$ and $r_2$ can be carried out for any regular pattern.  So that we may generalize the proof and get the following conclusion. 
~
\\ 
{\bf Conclusion 1.}
For arbitrary {\em regular} patterns $\mathbb{P}$ and $\mathbb{Q}$, and strictly monotonic function $f: C_{\mathbb{P}} \rightarrow C_{\mathbb{Q}}$:   
\begin{enumerate}
\item If $f^{\infty}/_{\sim}$ is a function, then $f$ is continuous. 
\item If $f$ is continuous, then $f^{\infty}/_{\sim}$ is also continuous in the classical sense, i.e. as the continuity of a function from  topological space $C^{\infty}_{\mathbb{P}}/_{\sim}$  into topological space $C^{\infty}_{\mathbb{Q}} /_{\sim}$.
\item It follows from (1) i (2) that if $f^{\infty}/_{\sim}$ is a function, then $f^{\infty}/_{\sim}$ is a continuous function. 
\end{enumerate}
The proof of (2) is simple and is as follows.  If $f^{\infty}/_{\sim}$ is not continuous, then there is sequence $(e_1,\ e_2,\ \dots e_n,\ \dots )$ converging to $e$ such that the sequence $(f^{\infty}/_{\sim}(e_1),\ f^{\infty}/_{\sim}(e_2),\ \dots f^{\infty}/_{\sim}(e_n),\ \dots )$ converges to $r \ \neq \ f^{\infty}/_{\sim}(e)$. 
Let $[u] = e$, and $[v] = f^{\infty}/_{\sim}(e)$, and $[u_n] = e_n$ for any $n$.  
Then, there are $m$ and $k$ such that for all $i>k$ and all $j > m$, the finite sequences $v(j)$ and $(f^{\infty}(u_i))(j)$ are not adjacent. Since $f$ is strictly monotonic, 
there is $l_j$ such that $v(j)$ is a prefix of $f(u(l_j))$, and there is $n_j$ such that $(f^{\infty}(u_i))(j)$ is a prefix of $f(u_i(n_j))$. By the regularity of the pattern $\mathbb{Q}$, also $f(u(l_j))$ and $f(u_i(n_j))$ are not adjacent. 
Since the sequence ($[u_i]$, $i= 1,2, \dots $) converges to  $[u]$\ \ (here the regularity of $\mathbb{P}$ is essential), for sufficiently large $i$, the finite sequences $u(l_j)$ and $u_i(n_j)$ are adjacent. Hence, function $f$ is not continuous, and the proof is completed. 

On the other hand, for any continuous function $\bar{f}$ on the unit interval $[0; 1]$ (in the classic sense), there is a strictly monotonic and continuous function $f: C_{\mathbb{E}^1} \rightarrow C_{\mathbb{E}^1}$, such that $f^{\infty}/_{\sim}$ and  
$\bar{f}$ are the same up to the homeomorphism between $C^{\infty}_{\mathbb{E}^1} /_{\sim}$ and the unit interval. 
The required function $f$ is defined as follows. 
For any object $x$ of $C_{\mathbb{E}^1}$, let $I_x$ denote the sub-interval of $[0; 1]$ corresponding to $x$.  
For any $x$, the let $f(x)$ be defined as $y$ such that $I_y$ is the smallest (with regard to inclusion) interval such that for any $r \in I_x$, \ \ $\bar{f}(r) \in I_y$. If such smallest interval does not exits, then by the continuity of $\bar{f}$, the function $\bar{f}$ is constant on $I_x$. Then, let $f(x)$ be defined as the first $y$ of length greater than the length of $x$, such that for any $r \in I_x$, \ \ $\bar{f}(r) \in I_y$. The function $f$ defined in this way is strictly  monotonic and continuous. The classic continuity of $\bar{f}$ is essential here. The extension of $f$, i.e $f^{\infty}/_{\sim}$ is the same (modulo the homomorphism) as $\bar{f}$. 

In the general case, the above statement is also true. That is, for arbitrary  {\em regular} patterns $\mathbb{P}$ and $\mathbb{Q}$, if $C^{\infty}_{\mathbb{P}} /_{\sim}$ corresponds (and is homeomorphic) to topological space $T_P$, and $C^{\infty}_{\mathbb{Q}} /_{\sim}$ corresponds (and is homeomorphic) to topological space $T_Q$, then for any continuous function 
$\bar{g}:T_P \rightarrow T_Q$ (in the classic sense), there is a strictly monotonic and continuous function $g: C_{\mathbb{P}} \rightarrow C_{\mathbb{Q}}$,  such that $g^{\infty}/_{\sim}$ and  
$\bar{g}$ are the same up to the homeomorphisms. \\ 
Proof is almost the same as for the Euclidean pattern $\mathbb{E}^1$. The only difference consists in definition of the set $I_x$. Now, it is defined as $[u] \in I_x$ if $x$ is a prefix of $u$. 

Note that the restriction of the function $g$ to the domain of objects (of type $C_{\mathbb{P}}$) of length $k$ may be seen as the $k$-th approximation of the continuous function $\bar{g}$.
%

The conclusion for this section is as follows. The continuity of strictly monotonic function $f$ is necessary and sufficient condition for  $f^{\infty}/_{\sim}$ to be a function, and a continuous function. 

\section{The Brouwer Continuity Theorem} 
\label{btc}

Let us recall that the unit interval $[0;\ 1]$ is homeomorphic to $C^\infty_{\mathbb{E}^1}/_{\sim}$. 

According to Brouwer, every real number (from the unit interval $[0; \ 1]$) is represented by a choice sequence that is an element of the set $C^\infty_{\mathbb{E}^1}$. 
However, there may be several different choice sequences (say $v$ and $u$) that correspond to the same real number $r$. Then, the sequences are adjacent (equivalent), i.e.  $v \sim u$, and determine one equivalence class (an element of the set $C^\infty_{\mathbb{E}^1}/_{\sim}$) that may be identified with the real number $r$.  

The Brouwer Continuity Theorem states that all functions on the unit interval are continuous, that is, any function that is defined everywhere on the unit interval of real numbers is also continuous. This is the famous ``paradox'' in intuitionistic Analysis. 
The original proof of the ``paradox theorem'' was based on the Brouwer's axiom of continuity. The axiom has been  presented in Section \ref{cs}. For the proof  see   \cite{veldman2001understanding}, and Stanford Encyclopedia of Philosophy \url{https://plato.stanford.edu/entries/intuitionism/}.  

Let us recall our Postulate, however, restricted only to pattern $\mathbb{E}^1$ that corresponds to the unit interval $[0; \ 1]$. 
\\
{\em Any function $F$ on choice sequences 
(i.e. from $C^\infty_{\mathbb{E}^1}$ 
into $C^\infty_{\mathbb{E}^1}$) is the extension of a strictly
monotonic function $f$ of type 
$C_{\mathbb{E}^1} \rightarrow C_{\mathbb{E}^1}$. That is, $F$ is one and the same as  $f^{\infty}$. }
\\
It means that in order to compute an initial segment (prefix) of the value $F(v)$, the function $F$ needs only an initial segment of its argument $v$. 

Since the unit interval $[0;\ 1]$ is homeomorphic to $C^\infty_{\mathbb{E}^1}/_{\sim}$, the Postulate implies that any function on the unit interval is of the form $f^{\infty}/_{\sim}$.   
For $f^{\infty}/_{\sim}$ to be a function, $f$ must be continuous. 
Form the Conclusion 1, the continuity of $f$ implies the continuity of $f^{\infty}/_{\sim}$ in the classic sense. 
\\
{\bf The Conclusion.} The Postulate implies that any function on the unit interval $[0;\ 1]$ is continuous.  

Note that (without assuming the Postulate) any classic continuous function  on the unit interval is (up the homeomorphism) the extension, and then the transformation of a strictly monotonic continuous function on $C_{\mathbb{E}^1}$; see the previous Section. 
Hence, the Postulate is satisfied by the continuous functions on the unit interval. However, the Postulate is for  {\em all} functions. 

We may conclude that The Brouwer Continuity Theorem is not a paradox if the real numbers (i.e. the interval $[0;\ 1]$) are grounded as the type $C_{\mathbb{E}^1}$, and functions on the unit interval are grounded as strictly monotonic and continuous functions of type $C_{\mathbb{E}^1} \rightarrow C_{\mathbb{E}^1}$.

In the original work of Brouwer (also unpublished manuscripts, see the PhD thesis by Johannes \cite{Kuiper} as a comprehensive source) on Continuum, there is no explicit adjacency relation. It seems that this very relation is the key to understand the Continuum and its grounding, as well as to explain the controversies of the intuitionistic Continuum  with the classic Mathematics.  

The Brouwer's continuity theorem can be generalized to functions on the topological spaces constructed according to arbitrary regular patterns. 

Note that  $f^{\infty}/_{\sim}$ is an abstract notion (its definition involves actual infinity) and its grounding (computational contents) is in the function $f$.
It seems that the strict monotonicity and continuity of $f$ is much more simple and primitive than the classic notion of continuity of $f^{\infty}/_{\sim}$. Moreover, function $f$ that is not continuous (and/or not strictly monotonic) may be also of interest. 

One may say that the proposed grounding is nothing but approximation of the Continuum, and a strictly monotonic and continuous function is an approximation of a continuous function on Continuum. If it is supposed that such Continuum {\em really does exist}, as well as the functions on that Continuum {\em do exist}, then indeed the proposed grounding provides only their approximations. 

Strictly monotonic and continuous functions on the types $C_{\mathbb{Q}}$, (where  $\mathbb{Q}$ is an arbitrary  regular pattern) are of particular interest in finite element method (numerical method for solving problems of engineering and mathematical physics), as well as in computer graphics where 2-dimensional grid of pixels, and 3-dimensional grid of voxels, correspond to $graph^{\mathbb{Q}}_k$, where $k$ determines the (2D and 3D) image resolution, finer if $k$ is bigger.  

Let us summarize the Section. 
The regular division and adjacency patterns, and the resulting types may be considered as a grounding for the Continuum -- 
{\em ``a concept in its own right and independent of number''}. 
This may be seen as a realization of the idea of Brouwer's Continuum.  

It seems that Professor Luitzen E. J. Brouwer was right and his intuitionism (when properly grasped and elaborated) constitutes the ultimate Foundations of Mathematics. 



\section{The Euclidean Continua} 

It seems that Indiscernibility condition, and the Homogeneity condition form Section \ref{Euclidean patterns} express the common intuitions related to the notion of Euclidean Continuum. 
So that, Euclidean Continuum is a generic concept independent of the classic notion of Euclidean spaces as the Cartesian product of a finite number of copies of the set of real numbers, i.e. $\mathbb{R}^n$. 

Recall that $\mathbb{E}^n$ denote the $n$-dimensional Euclidean pattern.  The type $C_{\mathbb{E}^n}$  is the grounding of the topological space $C^{\infty}_{\mathbb{E}^n} /_{\sim}$ that is isomorphic to the $n$-dimensional unit cube, i.e. $[0;\ 1]^n \subset \mathbb{R}^n$. 
Hence, it seems reasonable for the types resulting from the Euclidean patterns to be called {\em Euclidean Continua}. 

\begin{figure}[h]
	\centering
	\includegraphics[width=0.4\textwidth]{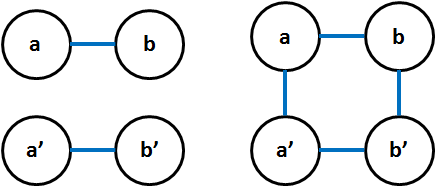}
	\caption{d-pattern for dimension 2 made from two d-patterns of dimension 1}
	\label{cartesianPattern}
\end{figure}

For $\mathbb{E}^2$, d-pattern is constructed by making a copy of the $1$ dimensional d-pattern and by joining to the original one in such a way that any original part and its copy become adjacent as shown in Fig. \ref{cartesianPattern}. 

Generally, for $\mathbb{E}^n$, make one copy of its d-pattern, and put it together. Any part and its copy are defined as adjacent.  
This corresponds to the Cartesian product. Give to copies different names. Then, the result is the d-pattern for $\mathbb{E}^{n+1}$. For the m-pattern, it is not so easy.  

The natural question arises: What is specific for the Euclidean patterns that they are the embedded base of human perception of space?  Note that the patterns are of dimension less than 4. For dimension 4 and higher our intuition fails. 

Note that using Euclidean patterns (of arbitrary dimension) a lot of interesting topological spaces can be constructed. For non-Euclidean patterns some of the topological spaces may be weired and surprising, whereas some spaces beautiful as fractals. 

\subsection{Conclusion}
\label{fin}

It was shown that the division patterns and resulting primitive types with adjacency relations constitute together the computational grounding of the Continuum. 
So that the relations seem to be the missed primitive concept in the type theory. 

If Robert Harper was right, i.e. \\ 
{\em ...
 type theory is and always has been a theory of computation on which the entire edifice of mathematics ought to be built.  ... ``},\\
 then the relations should become first class citizens in the type theory. 

The real numbers are the cornerstone of Calculus, Homotopy theory,
Riemannian geometry, and many others important subjects. The research on the grounding of Mathematics is continued in the paper {\em Asymptotic combinatorial constructions of Geometries} available at \url{https://arxiv.org/abs/1904.05173}. 

As the final conclusion, let me stress that there is nothing new in the approach presented above. All this (the so called intuitions) has been already known for ages in Mathematics and described formally
by definitions and theorems in the Analysis, Geometry, and Topology. The mere novelty of the proposed approach consists in explicit constructions of these very intuitions; the constructions that provide the grounding for these formal descriptions.

\bibliographystyle{ios2-nameyear}  
\bibliography{biblio-TO-ANG}     

\end{document}